\documentclass[12pt]{article}
\usepackage{amsfonts,amsthm,amsmath,amssymb}

\newtheorem{thm}{Theorem}[section]
\newtheorem{prop}[thm]{Proposition}
\newtheorem{lemma}[thm]{Lemma}
\newtheorem{cor}[thm]{Corollary}
\newtheorem{example}{Example}
\newtheorem{remark}[thm]{Remark}
\newtheorem{defin}{Definition}

\def\C{\mathbb{C}}
\def\R{\mathbb{R}}
\def\P{\mathbb{P}}

\def\g{\mathfrak{g}}
\def\h{\mathfrak{h}}
\def\gl{\mathfrak{gl}}
\def\h{\mathfrak{h}}
\def\c{\mathfrak{c}}

\def\V{V\oplus V}
\def\VV{V_1\oplus V_2}
\def\gs{\mathfrak{g}\oplus\mathfrak{g}}
\def\Ocl{\overline{\cal O}}
\def\DG{G\times G}

\def\d{\partial}

\def\p{\partial}
\def\t{\tilde}

\def\a{\alpha}

\def\phi{\varphi}
\def\s{\sigma}
\def\Ga{\Gamma}
\def\L{\Lambda}

\def\Oc{{\cal O}}

\def\S{{\cal S}}
\def\D{{\cal D}}

\def\dim{{\rm dim}}
\def\deg{{\rm deg}}

\def\Ad{{\rm Ad}}
\def\G{{\rm G}}

\def\T{{\rm T}}
\def\End{{\rm End}}

\def\Aut{{\rm Aut}}

\renewcommand\o{\overline}
\renewcommand\r{\rangle}
\renewcommand\l{\langle}

\title{Chern classes of reductive groups and an adjunction formula}
\author{Valentina Kiritchenko}
\date{}

\begin{document}
\maketitle
{\small In this paper, I construct noncompact analogs of the Chern
classes of equivariant vector bundles over complex reductive groups.
For the tangent bundle, these Chern classes yield an adjunction
formula for the Euler characteristic of complete
intersections in reductive groups. In the case where the complete intersection is a curve,
this formula gives an explicit answer for the Euler characteristic and the genus of the curve.}

\section{Introduction and main results}

Let $G$ be a connected complex reductive group. Consider a
faithful finite-dimensional representation $\pi: G\to GL(V)$ on
a complex vector space $V$. Let $H\subset\End(V)$ be a generic affine
hyperplane. The hypersurface $\pi^{-1}(\pi(G)\cap H) \subset G$ is
called a {\em hyperplane section} corresponding to the
representation $\pi$. The problem underlying this paper is how to
find the Euler characteristic of a hyperplane section or, more
generally, of the complete intersection of several hyperplane
sections corresponding to different representations.

The motivation to study such question comes from the case when the group $G=(\C^*)^n$ is
a complex torus. In this case, D.Bernstein,
A.Khovanskii and A.Kouchnirenko found an explicit and very beautiful answer in terms of
the weight polytopes of representations (see \cite{Khov}).
E.g. the Euler characteristic $\chi(\pi)$ of a hyperplane section corresponding to the representation $\pi$
is equal to $(-1)^n$ times the normalized volume of the weight polytope of $\pi$.
The proof uses an explicit relation between the Euler characteristic $\chi(\pi)$
and  the degree of the affine subvariety $\pi(G)$ in $\End(V)$:
$$\chi(\pi)=(-1)^{n-1}\deg~\pi(G).\eqno(1)$$
The degree is defined as usual. Namely, the degree of an affine subvariety $X\subset\C^N$ equals to
the number of the intersection points of $X$ with a generic affine subspace in $\C^N$ of complementary dimension.
For the degree $\deg~\pi(G)$ (that can also be interpreted as the self-intersection index of a hyperplane
section corresponding to the representation $\pi$) there is an explicit formula proved by
Kouchnirenko. Later D.Bernstein,
and Khovanskii found an analogous formula for the intersection index of hyperplane sections
corresponding to different representations.

How to extend these results to the case of arbitrary reductive
groups? It turned out that the formulas for the intersection
indices of several hyperplane sections can be generalized to
reductive groups and, more generally, to spherical homogeneous
spaces. For reductive groups, this was done by B.Kazarnovskii
\cite{Kaz}. Later, M.Brion established an analogous result for all
spherical homogeneous spaces \cite{Brion}. For reductive groups,
the Brion--Kazarnovskii theorem allows to compute explicitly the
intersection index of $n$ generic
hyperplane sections corresponding to  different representations.
The precise definition of the intersection index
is given in Section \ref{s.cp}.

However, when $G$ is an arbitrary reductive group, it is no longer true that $\chi(\pi)=(-1)^{n-1}\deg~\pi(G).$
K.Kaveh computed explicitly $\chi(\pi)$ and $\deg~\pi(G)$
for all representations $\pi$ of $SL_2(\C)$ . His computation shows
that, in general, there is a discrepancy between these two
numbers. Kaveh also listed some special representations of reductive groups,
for which these numbers still coincide
\cite{Kiu}.

In this paper, I will present a formula that, in particular,
generalizes formula (1) to the case of arbitrary reductive groups.
To do this I will construct algebraic subvarieties $S_i\subset G$, whose
degrees fill the gap between the Euler characteristic and the
degree. My construction is similar to one of the classical
constructions of the Chern classes of a vector bundle in the compact
setting (Subsection \ref{s.lbchern}). The subvarieties $S_i$ can be thought of as Chern classes
of the tangent bundle of $G$. I will also construct  Chern
classes of more general equivariant vector bundles over $G$
(Subsection \ref{s1.chern}). These Chern classes are in many
aspects similar to the usual Chern classes of compact manifolds.
There is an analog of the cohomology ring for $G$, where the Chern
classes of equivariant bundles live. This analog is the ring of
conditions constructed by C.De Concini and C.Procesi \cite{CP,C}(see
Section \ref{s.cp} for a reminder). It is useful in solving
enumerative problems. In particular, the intersection product in
this ring is well-defined.

I now formulate the main results. Denote by $n$ and $k$ the
dimension and the rank of $G$, respectively. Recall that the rank
is the dimension of a maximal torus in $G$. Denote by $[S_1],
\ldots, [S_{n}]$ the Chern classes of the tangent bundle of $G$ as
elements of the ring of conditions, and denote by $S_1$,\ldots,
$S_n$ subvarieties representing these classes. In the case of the
tangent bundle, it turns out (see Lemma \ref{dim}) that the
the higher
Chern classes $[S_{n-k+1}]$,\ldots, $[S_{n}]$ vanish. E.g. if $G$ is
a torus, then all Chern classes $[S_i]$ vanish.

Let $H_1$,\ldots, $H_m$ be a {\em generic} collection of $m$
hyperplane  sections corresponding to faithful representations
$\pi_1$,\ldots, $\pi_m$ of the group $G$ (for the precise meaning
of ``generic" see Subsection \ref{proofs}). Then the following
theorem holds.

\begin{thm}\label{compl}
The Euler characteristic of the complete intersection
$H_1\cap\ldots\cap H_m$ is equal to the term of degree $n$ in the
expansion of the following product:
$$(1+S_1+\ldots+S_{n-k})\cdot\prod_{i=1}^mH_i(1+H_i)^{-1}.$$
The product in this formula is the intersection product in the ring of conditions.
\end{thm}
This is very similar to the classical adjunction formula in
compact setting.

In particular, the Euler characteristic of just one hyperplane
section corresponding to a representation $\pi$ is equal to the
following alternating sum. Put $S_0=G$. Then
$$\chi(\pi)=\sum_{i=0}^{n-k}(-1)^{n-i-1}\deg~\pi(S_i).$$
The latter formula may have applications in the theory of
generalized hypergeometric equations. In the torus case,
I.Gelfand, M.Kapranov and A.Zelevinsky showed that the Euler
characteristic $\chi(\pi)$ gives the number of integral solutions
of the generalized hypergeometric system associated with the
representation $\pi$ \cite{GKZ}. A similar system can be
associated with the representation $\pi$ of any reductive group
\cite{Kapranov2}. In the reductive case, the number of integral
solutions of such a system is also likely to coincide with
$\chi(\pi)$.

The proof of Theorem \ref{compl}
is similar to the proof by Khovanskii \cite{Khov} in the
torus case.  Namely, Theorem \ref{compl} follows from the
adjunction formula applied to the closure of a complete
intersection in a suitable regular compactification of $G$ (see
Subsection \ref{proofs}). The key ingredient is a description of
the tangent bundles of regular compactifications due to Ehlers \cite{E} and Brion \cite{Brion4}. This description
is outlined in Subsection \ref{s.chern}.

The remaining problem is to describe the Chern classes
$[S_1],\ldots,[S_{n-k}]$ so that their intersection indices with
hyperplane sections may be computed explicitly. So far there is
such a description for the first and the last Chern classes (see
Subsection \ref{s.deg}). Namely, $[S_1]$ is the class of a generic
hyperplane section corresponding to the irreducible representation
with the highest weight $2\rho$. Here $\rho$ is the sum of all
fundamental weights of $G$. This description follows from a result of
A.Rittatore \cite{Rit} concerning the first Chern class of reductive group compactifications.
The last Chern class $[S_{n-k}]$ is up
to a scalar multiple the class of a maximal torus in $G$.
There is a hope that the intersection indices of other Chern classes $S_i$
with hyperplane sections can also be computed using a formula similar to the Brion--Kazarnovskii
formula.

If a complete intersection is a curve, i.e.
$m=n-1$, then the formula of Theorem \ref{compl} involves only the
first Chern class $[S_1]$. In this case, the computation of $[S_1]$ together
with the Brion--Kazarnovskii formula allows us to compute explicitly
the Euler characteristic and the genus of a curve in $G$ in terms of the weight polytopes of
$\pi_1$,..., $\pi_m$ (see Corollaries \ref{curve} and \ref{curve2}, Subsection \ref{proofs}).
Note that these two numbers completely describe the topological type of a curve.

Most of the constructions and results of this paper can be extended without any change to the case of
arbitrary spherical homogeneous spaces. This is discussed in Section \ref{s.sphere}.

I am very grateful to Mikhail Kapranov and Askold Khovanskii for numerous
stimulating discussions and suggestions. I would like to thank Kiumars Kaveh for useful
discussions and Michel Brion for valuable remarks on the first version of this paper.
I am also grateful to the referee for many useful remarks and comments.

Part of the results of this paper were included into my PhD thesis at
the University of Toronto \cite{K}.

Throughout this paper, whenever a group action is mentioned, it is always
assumed that a complex algebraic group acts on a complex algebraic variety by algebraic automorphisms.
In particular, by a homogeneous space for a group I will always mean
the quotient of the group by some closed algebraic subgroup.

The following remarks concern notations. In this paper, the term
equivariant (e.g. equivariant compactification, bundle, etc.) will
always mean equivariant under the action of the doubled group
$\DG$, unless otherwise stated.   The Lie algebra of $G$ is denoted
by $\g$. I also fix an embedding $G\subset GL(W)$ for some vector
space $W$. Then for $g\in G$ and $A\in\g$, notation $Ag$ and $gA$
mean the product of linear operators in $\End(W)$.

\section{Equivariant compactifications and the ring of
conditions} \label{s.cp} This section contains some well-known
notions and theorems, which will be used in the sequel. First, I
define the notion of spherical action and describe
equivariant compactifications of reductive groups following
\cite{CP2}, \cite{Kapranov2} and \cite{Tim}.
Then I state Kleiman's transversality theorem \cite{Kleiman} and
recall the definition of the ring of conditions \cite{CP,C}.

\paragraph{\bf Spherical action.} Reductive groups are partial cases of more general {\em spherical}
homogeneous spaces. They are defined as follows. Let $G$ be a connected
complex reductive group, and let $M$ be a homogeneous space under $G$.
The action
of $G$ on $M$ is called {\em spherical}, if a Borel subgroup of
$G$ has an open dense orbit in $M$. In this case, the homogeneous
space $M$ is also called spherical. An important and very useful
property, which characterizes a spherical homogeneous space $M$, is that any
compactification of $M$ equivariant under the action of $G$ contains only
a finite number of orbits \cite{Knop2}.

There is a natural action of the group $G\times G$ on $G$ by left
and right multiplications. Namely, an element $(g_1,g_2)\in
G\times G$ maps an element $g\in G$ to $g_1gg_2^{-1}$. This action
is spherical as follows from the Bruhat decomposition of $G$ with
respect to some Borel subgroup. Thus the group $G$ can be
considered as a spherical homogeneous space of the doubled group
$\DG$ with respect to this action. For any representation
$\pi:G\to GL(V)$ this action can be extended straightforwardly
to the action of $\pi(G)\times\pi(G)$ on the whole $\End(V)$ by
left and right multiplications. I will call such actions {\em
standard}.

\paragraph{\bf Equivariant compactifications.} With any representation $\pi$ one can associate
the following compactification of $\pi(G)$. Take the projectivization $\P(\pi(G))$ of $\pi(G)$ (i.e. the set
of all lines in $\End(V)$ passing through a point of $\pi(G)$ and the origin),
and then take its closure in $\P(\End(V))$.
We obtain a projective variety $X_\pi\subset\P(\End(V))$
with a natural action of $G\times G$ coming from the standard action of
$\pi(G)\times\pi(G)$ on $\End(V)$. Below I will list
some important properties of this variety.

Assume that $\P(\pi(G))$ is isomorphic to $G$. Fix a maximal torus $T\subset G$. Let $L_T$ be its
character lattice.
Consider all weights of the representation $\pi$, i.e. all characters
of the maximal torus $T$ occurring in $\pi$. Take their convex hull $P_\pi$ in $L_T\otimes\R$.
Then it is easy to see that $P_\pi$ is a polytope
invariant under the action of the Weyl group of $G$. It is called the {\em weight polytope}
of the representation $\pi$.
The polytope $P_\pi$ contains
information about the compactification $X_\pi$.
\begin{thm}\label{equiv}
1)  {\em (\cite{Tim}, Proposition 8)} The subvariety $X_\pi$ consists of a finite number of
$\DG$-orbits. These orbits are in one-to-one correspondence with
the orbits of the Weyl group acting on the faces of the polytope
$P_\pi$.

2)  Let $\sigma$ be another representation of $G$.
The normalizations of subvarieties $X_\pi$ and $X_\sigma$ are isomorphic if and only if
the normal fans corresponding to the polytopes $X_\pi$ and $X_\sigma$
coincide. If the first fan is a subdivision of the second, then there
exists an equivariant map from the normalization of $X_\pi$ to $X_\sigma$, and vice versa.
\end{thm}
The second part of Theorem \ref{equiv} follows from the general theory of spherical varieties
(see \cite{Knop2}, Theorem 5.1) combined with the description of compactifications $X_\pi$
via colored fans (see \cite{Tim}, Sections 7, 8).

In particular, suppose that the group $G$ is of adjoint type, i.e. the center of $G$ is trivial.
Let $\pi$ be an irreducible representation of $G$ with a strictly dominant highest weight.
It is proved in \cite{CP2} that the corresponding compactification $X_\pi$ of the group $G$
is always smooth and, hence, does not depend on the choice of a highest
weight. Indeed, the normal fan of
the weight polytope $P_\pi$ coincides with the fan of the Weyl chambers and their faces, so the second
part of Theorem \ref{equiv}
applies.
This compactification is
called the wonderful compactification and is denoted by $X_{can}$. It was introduced by De
Concini and Procesi \cite{CP2}. The  boundary divisor $X_{can}\setminus G$ is a divisor with normal crossings.
There are $k$ orbits $\Oc_1,\ldots,\Oc_k$ of codimension one in
$X_{can}$. The other orbits are obtained as the intersections of the
closures $\overline\Oc_1,\ldots,\overline\Oc_k$. More precisely,
to any subset $\{i_1,i_2,\ldots,i_m\}\subset\{1,\ldots,k\}$ there corresponds an
orbit
$\overline\Oc_{i_1}\cap\overline\Oc_{i_2}\cap\ldots\cap\overline\Oc_{i_m}$
of codimension $m$. So the number of orbits equals to $2^k$. There is a unique closed orbit
$\Ocl_1\cap\ldots\cap\Ocl_k$, which is isomorphic to the product of
two flag varieties $G/B\times G/B$. Here $B$ is a Borel subgroup
of $G$.

Compactifications of a reductive group arising from its
representations are examples of more general {\em equivariant
compactifications} of the group. A compact complex algebraic variety with an
action of $\DG$ is called an {\em equivariant
compactification} of $G$ if it satisfies the following conditions.
First, it contains an open dense orbit isomorphic to $G$. Second,
the action of $\DG$ on this open orbit coincides with the standard
action by left and right multiplications.

\paragraph{The ring of conditions.}
The following theorem gives a tool to define the intersection
index on a noncompact group, or more generally, on a homogeneous
space. Recall that two irreducible algebraic subvarieties $Y_1$ and $Y_2$ of an algebraic
variety $X$ are said to have {\em proper} intersection if either their intersection $Y_1\cap Y_2$ is empty
or all irreducible components of $Y_1\cap Y_2$ have dimension $\dim~Y_1+\dim~Y_2-\dim~X$.

\begin{thm}{\em (Kleiman's transversality theorem)} {\em \cite{Kleiman}}
\label{Kleiman} Let $H$ be a connected algebraic group, and let
$M$ be a homogeneous space under $H$. Take two algebraic subvarieties
$X,Y\subset M$. Denote by $gX$ the left translate of $X$ by an
element $g\in H$. There exists an open dense subset of $H$ such
that for all elements $g$ from this subset the intersection
$gX\cap Y$ is proper.
If $X$ and $Y$ are smooth, then $gX\cap Y$ is
transverse for general $g\in H$.

In particular, if $X$ and
$Y$ have complementary dimensions (but are not necessarily smooth),
then for almost all $g$ the translate $gX$ intersects $Y$ transversally at
a finite number of points, and this number does not depend on $g$.
\end{thm}
If $X$ and $Y$ have complementary dimensions, define the
{\em intersection index} $(X,Y)$ as the number $\#(gX\cap Y)$ of the intersection points
 for a generic $g\in H$. If one is interested in
solving enumerative problems, then it is natural to consider
algebraic subvarieties of $M$ up to the following equivalence. Two
subvarieties $X_1,X_2$ of the same dimension are equivalent if and
only if for any subvariety $Y$ of complementary dimension the
intersection indices $(X_1,Y)$ and $(X_2,Y)$ coincide. This
relation is similar to the numerical equivalence in algebraic
geometry (see \cite{Fulton2}, Chapter 19). Consider all formal linear combinations of algebraic
subvarieties of $M$ modulo this equivalence relation. Then the
resulting group $C^*(M)$ is called the {\em group of conditions} of $M$.

One can define an {\em intersection product} of two subvarieties $X,Y\subset M$ by
setting $X\cdot Y=gX\cap Y$, where $g\in G$ is generic. However,
the intersection product sometimes is not well-defined on the group of
conditions (see  \cite{CP} for a counterexample).
A remarkable fact is that for spherical homogeneous spaces the intersection product is
well-defined, i.e. if one takes different representatives of the
same classes, then the class of their product will be the same
\cite{CP,C}. The corresponding ring $C^*(M)$ is called the
{\em ring of conditions}.

In particular, the group of conditions $C^*(G)$ of a reductive
group is a ring. De Concini and Procesi related the
ring of conditions to the cohomology rings of equivariant
compactifications as follows. Consider the set $\S$ of
all smooth equivariant compactifications of the group $G$.
This set has a natural
partial order. Namely, a compactification $X_\s$ is greater
than $X_\pi$ if $X_\s$ {\em lies} over $X_\pi$,
i.e. if there exists a map $X_\s\to X_\pi$
commuting with the action of $\DG$. Clearly, such a map is  unique,
and it induces a map of cohomology rings $H^*(X_\pi)\to
H^*(X_\s)$.

\begin{thm}{\em \cite{CP,C}}\label{condi}
The ring of conditions $C^*(G)$ is isomorphic to the direct limit over the set $\S$
of the cohomology rings $H^*(X_\pi)$.
\end{thm}
De Concini and Procesi proved this theorem in \cite{CP} for
symmetric spaces. In \cite{C} De Concini noted that their
arguments go verbatim for arbitrary spherical homogeneous spaces.

\section{Chern classes of reductive groups}
\subsection{Preliminaries}\label{s.lbchern}

\paragraph{Reminder about the classical Chern classes.}
In this paragraph, I will recall one of the classical definitions
of the Chern classes, which I will use in the sequel. For more details
see \cite{GH}.

Let $M$ be a compact complex manifold, and let $E$ be a vector
bundle of rank $d$ over $M$. Consider $d$  global sections
$s_1,\ldots,s_d$ of $E$ that are $C^\infty$-smooth.
Define their {\em $i$-th degeneracy locus} as the set of all
points $x\in M$ such that the vectors $s_1(x),\ldots,s_{d-i+1}(x)$
are linearly dependent. The homology class of
the $i$-th degeneracy locus is the same for all {\em generic} choices of
the sections $s_1(x),\ldots, s_d(x)$ \cite{GH}. It is called the {\em $i$-th
Chern class} of $E$.

In what follows, I will only consider complex vector bundles that have plenty of algebraic global sections (so that
in the definition of the Chern classes, it will be possible to take only algebraic global sections instead of
$C^\infty$--smooth ones).

In particular, there
is the following way to choose generic global sections. Let $\Ga(E)$
be a finite-dimensional subspace in the space of all global
$C^\infty$-smooth sections of the bundle $E$. Suppose
that at each point $x\in M$ the sections of $\Ga(E)$ span the
fiber of $E$ at the point $x$. Then there is an open dense subset
$U$ in $\Ga(E)^d$ such that for any collection of global sections
$(s_1,\ldots,s_d)\subset U$ their $i$-th degeneracy locus is a representative
of the $i$-th Chern class of $E$.

I will also use the following classical construction that associates with
the subspace $\Ga(E)$ a map from the variety $M$ to a Grassmannian. Denote by $N$ the dimension
of $\Ga(E)$. Let $\G(N-d,N)$ be the Grassmannian of  subspaces of dimension $(N-d)$ in $\Ga(E)$.
One can map $M$ to $\G(N-d,N)$
by assigning to each point $x\in M$ the subspace of all sections
from $\Ga(E)$ that vanish at $x$.
By construction of the map the vector bundle $E$
coincides with the pull-back of the tautological quotient vector bundle
over the Grassmannian $\G(N-d,N)$. Recall that the tautological quotient vector
bundle over $\G(N-d,N)$ is the quotient of two bundles. The first one is the trivial vector bundle whose fibers
are isomorphic to $\Ga(E)$, and the second is the tautological vector bundle whose fiber at a point $\L\in\G(N-d,N)$
is isomorphic to the corresponding subspace $\L$ of dimension $N-d$ in $\Ga(E)$.

Using the definition of the Chern classes given above, it is easy to check
that the $i$-th Chern class of the tautological quotient vector
bundle is the homology class of the following Schubert cycle. Let
$\L^1\subset\ldots\subset\L^d\subset\Ga(E)$
be a partial flag such that $\dim~\L^j=j$. In the sequel, by a partial flag I will always mean
a partial flag of this type. The {\em $i$-th Schubert cycle} $C_i$ corresponding
to such a flag consists of all points $\L\in\G(N-d,N)$ such that the subspaces $\L$ and $\L^{d-i+1}$
have nonzero intersection.

The following proposition relates the Schubert cycles $C_i$
to the Chern classes of $E$.

\begin{prop}\label{lbchern}{\em \cite{GH}}
Let $p:M\to\G(N-d,N)$ be the map constructed above, and
let $C_i$ be the $i$-th Schubert cycle corresponding to a generic partial flag in $\Ga(E)$. Then the $i$-th
Chern class of $E$ coincides with the homology class of the inverse image of $C_i$ under the map $p$:
$$c_i(E)=[p^{-1}(C_i)].$$
\end{prop}
In particular, this proposition allows to relate the definition of the Chern classes
via degeneracy loci to other classical definitions.

In the sequel,  the following statement will be used. For any algebraic subvariety $X\subset\G(N-d,N)$,
a partial flag can be chosen in such a way that the corresponding Schubert cycle $C_i$ has
proper intersection with $X$. This follows from Kleiman's transversality theorem, since
the Grassmannian $\G(N-d,N)$ can be regarded as a homogeneous space under the natural action of the group
$GL_N$. Then any left translate of a Schubert cycle $C_i$ is again a Schubert cycle of the same type.

\paragraph{Equivariant vector bundles.}\label{s.eqbundle}
In this paragraph, I will recall the definition and some well-known properties
of equivariant vector bundles.

Let $E$ be a vector bundle of rank $d$ over $G$. Denote by
$V_g\subset E$ the fiber of $E$ lying over an element $g\in G$.
Assume that the standard action of $G\times G$ on $G$ can be
extended linearly to $E$.
More precisely, there exists a homomorphism $A:G\times
G\to\Aut(E)$ such that $ A(g_1,g_2)$ restricted to the
fiber $V_g$ is a linear operator from $V_g$ to $V_{g_1gg^{-1}_2}$.
If these conditions are satisfied, then the vector bundle $E$ is
said to be {\em equivariant} under the action of $G\times G$.

Two equivariant vector bundles $E_1$ and $E_2$ are equivalent if there
exists an isomorphism between $E_1$ and $E_2$ that is compatible
with the structure of fiber bundle and with the action of
$G\times G$. The following simple and well-known proposition describes
equivariant vector bundles on $G$ up to this equivalence relation.

\begin{prop}\label{bundle}
The classes of equivalent equivariant vector bundles of rank $d$
are in one-to-one correspondence with the linear representations
of $G$ of dimension $d$.
\end{prop}

Indeed, with each
representation $\pi:G\to V$ one can associate a bundle $E$ isomorphic to $G\times V$
with the following action of $\DG$: $$A(g_1,g_2):(g,v)\to(g_1gg^{-1}_2, \pi(g_1)v).$$
Then $A(g,g^{-1})$ stabilizes the identity element $e\in G$ and acts on the fiber $V_e=V$
by means of the operator $\pi(g)$.

E.g. the adjoint representation of $G$ on the Lie
algebra $\g=\T G_e$ corresponds to the tangent bundle $\T G$ on $G$.
This example will be important in the sequel.

Among all algebraic
global sections of an equivariant bundle $E$ there are two
distinguished subspaces, namely, the subspaces of left- and
right-invariant sections. They consists of sections that are
invariant under the action of the subgroups $G\times{e}$ and $e\times G$,
respectively. Both
spaces can be canonically identified with the vector space
$V$. Indeed, any vector $X\in V$ defines a
right-invariant section $v_r(g)=(g,X)$. Then it is easy to see that
any left-invariant section $v_l$ is given by the formula
$v_l(g)=(g,\pi(g)Y)$ for $Y\in V$.

Denote by $\Ga(E)$ the space of all global sections of
$E$ that are obtained as sums of left- and right-invariant
sections.
Let us find the dimension of the vector space $\Ga(E)$. Clearly, if
the representation $\pi$ does not contain any trivial
sub-representations, then $\Ga(E)$ is canonically isomorphic to
the direct sum of two copies of $V$. Otherwise, let $C\subset V$
be the maximal trivial sub-representation. Embed $C$ to $\V$
diagonally, i.e. $v\in C$ goes to $(v,v)$. It is easy to see that $\Ga(E)$ as a
$G$-module is isomorphic to the quotient space $(V\oplus V)/C$.
Denote by $c$ the dimension of $C$. Then the dimension of $\Ga(E)$
is equal to $2d-c$.

\subsection{Chern classes with values in the ring of conditions} \label{s1.chern}
In this subsection, I define Chern classes of equivariant vector
bundles over $G$. These Chern classes are elements of the ring of
conditions $C^*(G)$. Unlike the usual Chern classes in the compact
situation, they measure the complexity of the action of $G\times
G$ but not the topological complexity (topologically any $G\times
G$--equivariant vector bundle over $G$ is trivial). While the
definition of these classes does not use any compactification
it turns out that they are related to the usual
Chern classes of certain vector bundles over equivariant
compactifications of $G$.

Throughout this subsection, $E$ denotes the equivariant vector
bundle over $G$ of rank $d$ corresponding to a representation
$\pi:G\to GL(V)$. In the subsequent sections, I will only use the
Chern classes of the tangent bundle.

\paragraph{\bf Definition of the Chern classes.}
An equivariant vector bundle $E$ has a special class $\Ga(E)$ of algebraic global
sections. It consists of all global sections that can be represented as
sums of left- and right-invariant sections.

\begin{example}\em
If $E=\T G$ is the tangent bundle, then $\Ga(E)$ is a very natural
class of global sections. It consists of all vector fields
coming from the standard action of $\DG$ on $G$. Namely,
with any element $(X,Y)\in\gs$ one can associate a vector
field $v\in\Ga(E)$ as follows:
$$v(x)=\left.\frac{d}{dt}\right|_{t=0}[e^{tX}x e^{-tY}]=Xx-xY.$$
\end{example}
This example suggests that one represent elements of $\Ga(E)$ not
as sums but as differences of left- and right-invariant sections.

The space $\Ga(E)$ can be employed to define Chern classes of $E$
as usual. Take $d$ generic sections $v_1,\ldots,v_d\in\Gamma(E)$.
Then the $i$-th Chern class is the $i$-th degeneracy locus of
these sections. More precisely, the $i$-th Chern class
$S_i(E)\subset G$ consists of all points $g\in G$ such that the first $d-i+1$
sections $v_1(g),\ldots,v_{d-i+1}(g)$ taken at $g$ are linearly dependent.
This definition almost repeats
one of the classical definitions of the Chern classes in the compact setting (see
Subsection \ref{s.lbchern}). The only difference is that global sections used in
this definition are not generic in the space of all sections. They are
generic sections of the special subspace $\Ga(E)$. If one drops this restriction and applies
the same definition, then the result will be trivial, since the bundle $E$ is topologically
trivial. In some sense, the Chern classes will sit at infinity in this case (the precise meaning will become clear
from the second part of this subsection). The purpose of my definition is to pull them back
to the finite part.

Thus for each $i=1,\ldots,d$ we get a family $\S_i(E)$ of
algebraic subvarieties $S_i(E)$ parameterized by collections of $d-i+1$ elements
from $\Ga(E)$. In the compact situation, all generic members of
an analogous family represent the same class in the cohomology
ring. The same is true here, if one uses the ring of conditions as
an analog of the cohomology ring in the noncompact setting.
\begin{lemma}\label{class}
For all collections $v_1,\ldots,v_{d-i+1}$ belonging to some open dense
subset of $(\Ga(E))^{d-i+1}$ the class of the corresponding
subvariety $S_i(E)$ in the ring of conditions $C^*(G)$ is the
same.
\end{lemma}

The lemma implies that the family $\S_i(E)$ parameterized by
elements of $(\Ga(E))^{d-i+1}$ provides a well-defined class
$[S_i(E)]$ in the ring of conditions $C(G)$.

\begin{defin}\label{chern}
The class $[S_i(E)]\in C^*(G)$ defined by the family
 $\S_i(E)$ is called the $i$-th Chern class of a vector
bundle $E$ with value in the ring of conditions.
\end{defin}

Before proving the lemma let me give another description of the Chern
classes $[S_i(E)]$.

\paragraph{\bf Maps to Grassmannians.}
In this paragraph, I apply the classical construction discussed in
Subsection \ref{s.lbchern} to define a map from the group $G$ to
the Grassmannian $G(d-c,\Ga(E))$ of  subspaces of dimension
$(d-c)$ in the space $\Ga(E)$.  Recall that $c$ is the dimension
of the maximal trivial sub-representation of $V$, and the
dimension of $\Ga(E)$ is $2d-c$ (see the end of Subsection
\ref{s.lbchern}).

Note that the global sections from the subspace $\Ga(E)$ span the
fiber of $E$ at each point of $G$. Hence, one can define a map
$\phi_E$ from $G$ to the Grassmannian $G(d-c,\Ga(E))$ as follows.
A point $g\in G$ gets mapped to the subspace $\L_g\subset\Ga(E)$
spanned by all global sections that vanish at $g$. Clearly, the
dimension of $\L_g$ equals to $(\dim~\Ga(E)-d)=(d-c)$ for all
$g\in G$. We get the map
$$\phi_E:G\to \G(d-c,\Ga(E));\quad \phi_E:g\mapsto\L_g.$$

The subspace $\L_g$ can be alternatively described using the graph
of the operator $\pi(g)$ in $\V$. Namely, it is easy to check that
$\L_g=\{(X,\pi(g)X), X\in V\}/C$. Then $\phi_E$ comes from the
natural map assigning to the operator $\pi(g)$ on $V$ its graph in
$\V$.

Clearly, the pull-back of the tautological quotient vector bundle
over $\G(d,\Ga(E))$ is isomorphic to $E$. Hence, the Chern
class $S_i(E)$ constructed via elements $v_1,\ldots,v_d$ is the inverse image of the
Schubert cycle $C_i$ corresponding to the partial flag $\l
v_1\r\subset\l v_1,v_2\r\subset\ldots\subset\l
v_1,\ldots,v_d\r\subset\Ga(E)$ (see Subsection \ref{s.lbchern}).
Here $\l v_1,\ldots,v_i\r$ denotes
the subspace of $\Ga(E)$ spanned by the vectors $v_1,\ldots,v_i$.

\begin{remark} \label{second}\em
This gives the following equivalent definition of $S_i(E)$. The Chern class $S_i(E)$
consists of all elements $g\in
G$ such that the graph of the operator $\pi(g)$ in $\V$ has a
nontrivial intersection with a generic subspace of dimension $d-i+1$
in $\V$.
\end{remark}

In particular, if the representation
$\pi:G\to GL(V)$ corresponding to a vector bundle $E$ has a
nontrivial kernel, then the $S_i(E)$ are invariant under left and
right multiplications by the elements of the kernel (since this is already true
for the preimage $\phi_E^{-1}(\L)$ of any point $\L\in\phi_E(G)$). E.g. the
Chern classes $S_i(\T G)$ are invariant under  multiplication
by the elements of the center of $G$.

We can now relate the
Chern classes $S_i(E)$ to the usual Chern classes of a vector
bundle over a compact variety.

Denote by $X_E$ the closure of $\phi_E(G)$ in the Grassmannian
$G(d-c,\Ga(E))$, and denote by $E_X$ the restriction of
the tautological quotient vector bundle to $X_E$.
We get a vector bundle on a compact variety. The $i$-th Chern class
of $E_X$ is the homology class of $C_i\cap X_E$
for a generic Schubert cycle $C_i$ (see Proposition \ref{lbchern}).
By Kleiman's transversality theorem applied to the Grassmannian $\G(d-c,\Ga(E))$
(see Subsection \ref{lbchern}), a generic Schubert
cycle $C_i$ has a proper intersection with
the boundary divisor $X_E\setminus\phi_E(G)$.
Hence, there is the following relation between the Chern classes of $E_X$
and generic members of the family $\S_i(E)$.

\begin{prop}\label{relation}
For a generic $S_i(E)$ the homology class of the closure
of $\phi_E(S_i(E))$ in $X_E$
coincides with the $i$-th Chern class of $E_X$.
\end{prop}
Thus the Chern classes $[S_i(E)]$ can be described via the usual Chern classes of the bundle $E_X$ over
the compactification $X_E$.

Let us study the variety $X_E$ in more detail. It
is a $\DG$--equivariant compactification of the group $\phi_E(G)$.
Indeed, the action of $\DG$ on $\phi_E(G)$ can be extended to the
Grassmannian $\G(d,\Ga(E))$ as follows. Identify $\Ga(E)$ with
$(\V)/C$ (see the end of Subsection \ref{s.lbchern}). The doubled group $\DG$ acts on $\V$ by means of the
representation $\pi\oplus\pi$, i.e. $(g_1,g_2)(v_1,v
_2)=(g_1v_1,g_2v_2)$ for $g_1,g_2\in G,v_1,v_2\in V$. The subspace $C\subset\V$ is invariant under this action.
Hence, the group $\DG$ acts on $\Ga(E)$. This action
provides an action of $\DG$ on the Grassmannian $\G(d-c,\Ga(E))$.
Clearly, the subvariety $X_E$ is invariant under this action.

{\bf Example 1 (Demazure embedding).} Let $G$ be a group of
adjoint type, and let $\pi$ be its adjoint representation on the
Lie algebra $\g$. The corresponding vector
bundle $E$ coincides with the tangent bundle of $G$. The
corresponding map $\phi_E:G\to\G(n,\gs)$ coincides with the
embedding constructed by Demazure \cite{CP2}. The Demazure map takes
an element $g\in G$ to the Lie subalgebra
$\g_g=\{(gXg^{-1},X),
X\in\g\}\subset\gs$. Clearly, the Demazure map provides an
embedding of $G$ into $\G(n,\gs)$.

It is easy to check that the Lie subalgebra $\g_g$ is the Lie algebra
of the stabilizer of an element $g\in G$ under the standard action of
$\DG$. Thus for any $A\in\g_g$ the corresponding vector field
vanishes at $g$, and the Demazure embedding coincides with $\phi_E$.
The compactification $X_E$ in this case is
isomorphic to the wonderful compactification $X_{can}$ of the group $G$ \cite{CP2}.
In particular, it is smooth.

\begin{defin}\label{Dem}
Let $G$ and $E$ be  as in Example 1. The restriction
of the tautological quotient vector bundle to $X_E\simeq X_{can}$ is
called the Demazure bundle and is denoted by $V_{can}$.
\end{defin}

If $E$ is the tangent bundle, then Proposition \ref{relation}
implies that the Chern class $S_i(E)$ is the inverse image of the
usual $i$-th Chern class of the Demazure bundle. The Demazure bundle
is considered in \cite{Brion4}, where it is related to the tangent bundles
of regular compactifications of the group $G$.

\begin{example}\em
{\bf a)} Let $G$ be $GL(V)$ and let $\pi$ be its tautological
representation on the space $V$ of dimension $d$. Then $\phi_E$ is
an embedding of $GL(V)$ into the Grassmannian $\G(d,2d)$. Notice
that the dimensions of both varieties are the same. Hence, the
compactification $X_E$ coincides with $\G(d,2d)$.

{\bf b)} Take $SL(V)$ instead of $GL(V)$ in the previous example.
Its compactification $X_E$ is a hypersurface in the Grassmannian
$\G(d,2d)$ which can be described as a hyperplane section of the
Grassmannian in the Pl\" ucker embedding. Consider the Pl\" ucker
embedding $p:\G(d,2d)\to\P(\L^d(\VV))$, where $V_1$ and $V_2$ are two
copies of $V$. Then $p(X_E)$ is a special hyperplane section of
$p(\G(d,2d))$. Namely, the decomposition $\VV$ yields a
decomposition of $\L^d(\VV)$ into a direct sum. This sum
contains two one-dimensional components $p(V_1)$ and $p(V_2)$
(which are considered as lines in $\L^d(\VV)$). In particular, for
any vector in $\L^d(\VV)$ it makes sense to speak of its
projections to $p(V_1)$ and $p(V_2)$.  On $V_1$ and $V_2$ there
are two special $n$-forms, preserved by $SL(V)$. These forms give
rise to two $1$-forms $l_1$ and $l_2$ on $p(V_1)$ and $p(V_2)$,
respectively. Consider the hyperplane $H$ in $\L^d(\VV)$ consisting
of all vectors $v$ such that the functionals $l_1$ and $l_2$ take the
same values on the projections of $v$ to $p(V_1)$ and $p(V_2)$,
respectively. Then it is easy to check that
$p(X_E)=p(\G(d,2d))\cap \P(H)$.
\end{example}

In the next section, I will be concerned with the case when $E=\T
G$ is the tangent bundle. In this case, the vector bundle $E_X$ is
closely related to the tangent bundles of regular
compactifications of the group $G$. Let us discuss this case in
more detail.

\begin{example}\em \label{Dem2}
This example is a slightly more general version of Example
1.
Let $\g=\g'\oplus\c$ be the decomposition of the Lie
algebra $\g$ into the direct sum of the semisimple and the central
subalgebras, respectively. Denote by $c$ the dimension of the
center $\c$. Let $E=\T G$ be the tangent bundle on $G$.
Then $\phi_E$ maps $G$ to the Grassmannian
$G(n-c,(\gs)/\c)$. It is easy to show that the image
of the map $\phi_E$ coincides with the adjoint group of $G$ and the image contains only
subspaces that belong to $(\g'\oplus\g')\subset(\gs)/\c$. Comparing this with Example 1,
one can easily see that
$X_E$ is
isomorphic to the wonderful compactification $X_{can}$ of the
adjoint group of $G$.

In this case, the bundle $E_X$ is the direct sum
of the Demazure bundle and the trivial vector bundle of rank $c$.
Indeed, for any subspace $\L_x\in X_E\simeq X_{can}\subset\G(n-c,\Ga(E))$ its
intersection with the subspace $\c^-=\{(c,-c),
c\in\c\}\subset\Ga(E)$ is trivial. Hence, the quotient space $\Ga(E)/\L_x$ coincides
with the direct sum $((\g'\oplus\g')/\L_x)\oplus\c^-$.

\end{example}

\paragraph{\bf Proof of Lemma \ref{class}}
The proof of Lemma \ref{class} relies on the following
fact. Let $Y_1$ and $Y_2$ be two subvarieties of codimension $i$
in the group $G$. Using Kleiman's transversality theorem and
continuity arguments, it is easy to show  that $Y_1$ and $Y_2$
represent the same class in the ring of conditions $C^*(G)$ if
there exists an equivariant compactification $X$ of the group $G$
such that the closures of $Y_1, Y_2$ in $X$ have  proper intersections
with all $\DG$--orbits (see \cite{CP} for the proof).

In particular, to prove  Lemma \ref{class} it is enough to produce
an equivariant compactification $X$ such that the closure of a
generic $S_i(E)$ has proper intersections with all $\DG$--orbits
in $X$. I claim that the compactification $X_E$ discussed in the previous paragraph
(see Proposition \ref{relation}) satisfies this condition.

Indeed, the closure of any $S_i(E)$ in $X_E$
coincides
with the intersection of $X_E$ with the Schubert cycle $C_i$ corresponding to a partial flag
in $\Ga(E)$. By Kleiman's transversality theorem applied to the Grassmannian $\G(d-c,\Ga(E))$
(see Subsection \ref{lbchern}), a partial flag can be chosen in such a way that the corresponding Schubert cycle has
proper intersections with all $\DG$--orbits in $X_E$. All partial flags with such property form an open
dense subset in the space of all partial flags.
Hence, for generic flags the corresponding
subvarieties $S_i$ represent the same class in the ring of conditions.

In the sequel, $S_i(E)$ will denote any subvariety of the
family $\S_i(E)$ whose class in the ring of conditions coincides
with the Chern class $[S_i(E)]$.

{\bf Remark.}
Recall that the ring of conditions $C^*(G)$ can be identified with
the direct limit of cohomology rings of equivariant
compactifications of $G$ (see Theorem \ref{condi}). It follows that under this
identification the Chern class $[S_i(E)]\in C^*(G)$ corresponds to
an element in the cohomology ring of the compactification $X_E$. In
particular for an adjoint group $G$, the Chern class $[S_i(\T X)]$
of the tangent bundle corresponds to some cohomology class of the
wonderful compactification of $G$.

\paragraph{ \bf Properties of the Chern classes of reductive groups.}
 The next lemma computes the dimensions of the
Chern classes. It also shows that if $G$ acts on $V$
without an open dense orbit, then the higher Chern classes automatically
vanish.

For any representation $\pi:G\to GL(V)$, there exists an open dense $G$--invariant subset in $V$ such that
the stabilizers of any two elements from this subset are conjugate subgroups of $G$  (see \cite{Rich}). In particular,
all elements from this subset
have isomorphic $G$--orbits. Such orbits are called {\em principal}.
Denote by $d(\pi)$ the dimension of a principal orbit of
$G$ in $V$.  If $G$ has an open dense
orbit in $V$, then $d(\pi)=d$. In my main example, when $\pi$ is the
adjoint representation, $d(\pi)=n-k$.

\begin{lemma}\label{dim} If $i>d(\pi)$, then $S_i(E)$ is empty, and
if $i\le d(\pi)$ then the dimension of $S_i(E)$ is equal to $n-i$.
\end{lemma}

\begin{proof}
Recall that $S_i(E)$ is the inverse image of $C_i$ under the map
$\phi_E:G\to \G(d-c,\Ga(E))$. Here $C_i$ is the $i$-th Schubert cycle corresponding
to a generic partial flag in $\Ga(E)$.
The codimension of $C_i$ in the Grassmannian $\G(d-c,\Ga(E))$ is equal to $i$.
Hence, by Kleiman's transversality theorem applied to  $\G(d-c,\Ga(E))$
(see Subsection \ref{lbchern}), the intersection $C_i\cap\phi_E(G)$ is either empty or proper and
has codimension $i$ in $\phi_E(G)$.
Then $S_i(E)=\phi_E^{-1}(C_i\cap\phi_E(G))$ is either empty or has codimension $i$ in $G$, because all fibers of the map $\phi_E$
are isomorphic to each other (each of them is isomorphic to the kernel of $\pi$). It remains to find out
all $i$ for which $S_i(E)$ is empty.

By Remark \ref{second}, the Chern class
$S_i(E)$ consists of all elements $g\in G$
such that the graph $\Ga_g=\{(v,\pi(g)v),v\in V\}\subset\V$ of
$\pi(g)$ has a nontrivial intersection with a generic subspace
$\L^{d-i+1}$ of dimension $d-i+1$ in $\V$. For all
$g\in S_i(E)\setminus S_{i+1}(E)$ the intersection $\Ga_g\cap\L^{d-i+1}$ has
dimension 1. Indeed, if $\dim(\Ga_g\cap\L^{d-i+1})\ge2$, then $\dim(\Ga_g\cap\L^{d-i})\ge1$
(since the subspace $\L^{d-i}\subset\L^{d-i+1}$ has codimension one in $\L^{d-i+1}$), and $g$ belongs to
$S_{i+1}(E)$.
Hence, there is a well-defined map
$$p:S_i(E)\setminus S_{i+1}(E)\to \P(D\cap \L^{d-i+1});\quad p:g\mapsto\P(\Ga_g\cap\L^{d-i+1}).$$
Here $D\subset\V$ is the union of all graphs $\Ga_g$ for $g\in
G$. In particular, the Chern class $S_i(E)$ is nonempty if and only if
$\P(D\cap\L^{d-i+1})$ is nonempty.

We now estimate the dimension of $D\cap\L^{d-i+1}$.
Since $D$ is not a variety, it is more convenient to take its Zariski closure $\o D$.
The subvariety $\o D$ is the closure of the image of the following morphism:
$$F: G\times V\to V\times V; \quad F: (g,v)\mapsto (v,\pi(g)v).$$
The source space $G\times V$ is an irreducible variety of dimension $n+d$, and
the general fibers of $F$ are isomorphic to the principal stabilizers, of dimension $n-d(\pi)$.
Hence $\dim~\o D=d+d(\pi)$, that is, $\o D$ has codimension $d-d(\pi)$.

Next, observe that $D$ is a constructible set, invariant under scalar multiplication. Hence
it contains a dense open subset (also invariant under scalar multiplication)
of the irreducible variety $\o D$.
Thus a general vector space $\L^{d-i+1}$ satisfies
$\dim(\o D\cap \L^{d-i+1})=d(\pi)-i+1$, if $i\le d(\pi)$,
and $D\cap \L^{d-i+1}$ is dense in this intersection.
In particular, if $i=d(\pi)$, then
$D\cap \L^{d-i+1}$ consists of several lines whose number is equal
to the degree of $\o D$. If $i>d(\pi)$, then $\o D\cap \L^{d-i+1}$
contains only the origin. It follows that if $i>d(\pi)$, then
$S_i(E)$ is empty.

\end{proof}

This proof also implies the following corollary. Denote by
$H\subset G$ the stabilizer of an element in a principal orbit of $G$ in $V$. The subgroup $H$ is
defined up to conjugation so its class in the ring of conditions is well-defined.
\begin{cor}\label{union}
An open dense subset of the subvariety $S_i(E)$ admits {\em almost} a fibration whose fibers are
translates of $H$. Here {\em almost} means that the intersection of different fibers always
lies in $S_{i+1}(E)\subset S_i(E)$.
In particular,
the last Chern class $S_{d(\pi)}(E)$ admits a true fibration and coincides with the disjoint
union of several translates of $H$. Their number equals to  the degree of a generic principal orbit of $G$
in $V$.
\end{cor}
The last statement follows from the fact that the degree of $D$ in $\V$
(see the proof of Lemma \ref{dim}) is equal to the degree of a generic principal orbit of $G$ in $V$.

In particular, let $E$ be the tangent bundle. Then the stabilizer
of a generic element in $\g$ is a maximal torus in $G$. Hence,
the last Chern class $S_{n-k}(\T G)$ is the
union of several translates of a maximal torus. The number of translates
is the cardinality of the Weyl group (the degree of a general orbit in the adjoint representation).

\subsection{The first and the last Chern classes}\label{s.deg}
Throughout the rest of the paper, I will only consider the Chern
classes $S_i=S_i(\T G)$ of the tangent bundle unless otherwise
stated. Theorem \ref{compl} expresses the Euler characteristic of
a complete intersection via the intersection indices of the Chern
classes $S_i$ with generic hyperplane sections. The  question
is how to compute these indices. If $[S_i]$ is a linear combination
of complete intersections of {\em generic} hyperplane sections corresponding
to some representations of $G$, then the answer to this question is given
by the Brion--Kazarnovskii formula. A hyperplane
section corresponding to the representation $\pi$  is called {\em generic}
if its closure in the compactification $X_\pi$ has proper intersections with all
$\DG$--orbits in $X_\pi$.

In this subsection, I
describe $S_1$ as a generic hyperplane
section. The description follows  from a result of
Rittatore \cite{Rit}. One can also compute the intersection
indices with the last Chern class $S_{n-k}$,
because $S_{n-k}$ is the union of translates of a maximal torus (see
Corollary \ref{union}). However, it seems that in
general the Chern class $S_i$, for $i\ne1$, is not a sum of complete
intersections. E.g. I can show that for $G=SL_3(\C)$ the Chern
class $[S_3]$ does not lie in the subring of $C^*(G)$ generated
by the classes of hypersurfaces.

\paragraph{Description of $S_1$.} The result of
 Rittatore for the first Chern class of regular compactifications (see \cite{Rit}, Proposition
4) implies that the class $[S_1]$ in
the ring of conditions can be represented by the doubled sum of the closures of all codimension one
Bruhat cells in $G$. Below I will deduce this description directly from the definition of $S_1$.

It is easy to show that $S_1\subset G$ is given by the equation
${\rm det}(\Ad(g)-A)=0$ for a generic $A\in\End(\g)$. Indeed,
the first Chern class $S_1(E)$ of any
equivariant vector bundle $E$ over $G$ consists of elements $g\in
G$ such that the graph of the operator $\pi(g)$ in $\V$ has a
nontrivial intersection with a generic subspace of dimension $n$
in $\V$ (see Remark \ref{second}). As a generic subspace, one can take the graph of a generic
operator $A$ on $V$. Then the graphs of operators $\pi(g)$ and $A$
have a nonzero intersection if and only if the kernel of the
operator $\pi(g)-A$ is nonzero.

The function ${\rm det}(\Ad(g)-A)$ is a linear combination of
matrix coefficients corresponding to all exterior powers of the adjoint
representation. Hence, the equation of $S_1$ is the equation of a
hyperplane section corresponding to the sum of all exterior powers
of the adjoint representation. Denote this representation by
$\sigma$. It is easy to check that the weight polytope $P_\sigma$
coincides with the weight polytope of the irreducible
representation $\theta$ with the highest weight $2\rho$ (here
$\rho$ is the half sum of all positive roots, or equivalently the sum of all fundamental weights). It
remains to prove that $S_1$ is generic, which means that the
closure of $S_1$ in $X_\sigma$ intersects all $\DG$--orbits along
subvarieties of codimension one. The proof of Lemma
\ref{class} implies that this is true for the wonderful
compactification, and the normalization of $X_\sigma$ is the wonderful
compactification by Theorem \ref{equiv} (since $P_\theta=P_\sigma$).

It is now easy to show that the doubled sum of the
closures of all codimension one
Bruhat cells in $G$ is equivalent to $S_1$. This is because the closures of  codimension one Bruhat cells
are generic hyperplane sections corresponding to the irreducible representations with fundamental
highest weights.

\paragraph{Description of $S_{n-k}$.}

By Corollary \ref{union} the last Chern class $S_{n-k}$ is the
disjoint union of translates of a maximal torus. Their number is equal to the
degree of a generic adjoint orbit in $\g$. The latter is equal to
the order of the Weyl group $W$. Denote by $[T]$ the class of a
maximal torus in the ring of conditions $C^*(G)$. Then the
following identity holds in $C^*(G)$:
$$[S_{n-k}]=|W|[T].$$
The degree of $\pi(T)$ can be computed using the formula of D.Bernstein,
Khovanskii and Koushnirenko \cite{Khov}.

\subsection{Examples}\label{s.exs}
\paragraph{$\bf G=SL_2(\C)$.} Consider the tautological embedding of $G$, namely,
$G=\{(a,b,c,d)\in\C^4:ad-bc=1\}$. Since the dimension of $G$ is 3
and the rank is 1,  by Lemma \ref{dim} we get that there are
only two nontrivial Chern classes: $S_1$ and $S_2$. Let us apply
the results of the preceding subsection to find them. The first
Chern class $S_1$ is a generic hyperplane section corresponding to
the second symmetric power of the tautological representation,
i.e. to the representation $\theta:SL_2(\C)\to SO_3(\C)$. In other
words, it is the intersection of $SL_2(\C)$ with a generic quadric in
$\C^4$. The second Chern class $S_2$ (which is also the last one
in this case) is the union of two translates of a maximal torus (or the intersection
of $S_1$ with a hyperplane in $\C^4$).

Let $\pi$ be a faithful representation of $SL_2(\C)$. It is a direct
sum of irreducible representations. Any irreducible
representation of $SL_2(\C)$ is isomorphic to the $i$-th symmetric
power of the tautological representation for some $i$.  Its weight
polytope is the line segment $[-i,i]$. Hence the weight polytope of
$\pi$ is the line segment $[-n,n]$ where $n$ is the greatest
exponent of symmetric powers occurring in $\pi$. Then the matrix coefficients
of $\pi$ are polynomials in $a,b,c,d$ of degree $n$. In
this case, it is easy to compute the degrees of subvarieties
$\pi(G)$, $\pi(S_1)$ and $\pi(S_2)$ by the Bezout theorem. Then
$\deg~\pi(G)=2n^3$, $\deg~\pi(S_1)=4n^2$, $\deg~\pi(S_2)=4n$.
Also, if one takes another faithful representation $\s$ with the
weight polytope $[-m,m]$, then the intersection index
 of $S_1$ with two generic hyperplane sections
corresponding to $\pi$ and $\s$, equals to $4mn$.

Since by Theorem \ref{compl} the Euler characteristic $\chi(\pi)$
of a generic hyperplane section is equal to
$\deg~\pi(G)-\deg~\pi(S_1)+\deg~\pi(S_2)$, we get
$$\chi(\pi)=2n^3-4n^2+4n.$$
This answer was first obtained by Kaveh who used
different methods \cite{Kiu}.

If $\pi$ is not faithful, i.e. $\pi(SL_2(\C))=SO_3(\C)$, consider
$\pi$ as a representation of $SO_3(\C)$. Then $\chi(\pi)$ is two
times smaller and equals to $n^3-2n^2+2n$.

Apply Theorem \ref{compl} to a curve $C$ that is the complete
intersection of two generic hyperplane sections corresponding to
the representations $\pi$ and $\s$. Then
$$\chi(C)=H_\pi\cdot H_\s\cdot H_\theta-H_\pi\cdot H_\s\cdot (H_\pi+H_\s)=-2mn(m+n-2).$$

\paragraph{$\bf G=(\C^*)^n$ is a complex torus.}
In this case, all left-invariant vector fields are also
right-invariant since the group is commutative. Hence, they are
linearly independent at any point of $G=(\C^*)^n$ as long as their
values at the identity are linearly independent. It follows that
all subvarieties $S_i$ are empty, and all the Chern classes
vanish. Then Theorem 1 coincides with a theorem of D.Bernstein
and Khovanskii \cite{Khov}.

\section{Chern classes of regular compactifications and proof of Theorem \ref{compl}}
\label{s.chernclass}
\subsection{Preliminaries}\label{s.prelim}
\paragraph{Chern classes of the tangent bundle.}
In this paragraph, I explain a method from \cite{E}, which in some
cases allows to find the Chern classes of smooth varieties.

Let $X$ be a smooth complex variety of dimension $n$, and let
$D\subset X$ be a divisor.
Suppose that $D$ is the union of $l$
smooth irreducible hypersurfaces $D_1,\ldots,D_l$ with normal crossings.
One can relate the tangent bundle $\T X$ of $X$ to the logarithmic tangent bundle,
consisting of those vector fields that preserve the divisor $D$.

Let
$L_X(D_1),\ldots,L_X(D_l)$ be the line bundles over $X$
associated with the hypersurfaces $D_1,\ldots,D_l$, respectively. I.e. the first
Chern class of the bundle $L_X(D_i)$ is the homology class of $D_i$.
One can also associate with $D$ the logarithmic tangent bundle
$V_X(D)$.  It is a holomorphic vector bundle over $X$ of rank $n$
that is uniquely defined by the following property. The holomorphic sections of
$V_X(D)$
over an open subset $U\subset X$ consist of all holomorphic vector fields
$v(x)$ on $U$ such that $v(x)$ restricted to $U\cap D_i$ is tangent to the hypersurface
$D_i$ for any $i$. The precise definition is as follows. Cover $X$ by local charts. If a chart
intersects the divisors $D_{i_1},\ldots, D_{i_k}$ choose local coordinates $x_1,\ldots,x_n$ such that
the equation of $D_{i_j}$ in these coordinates is  $x_j=0$. Then $V_X$ is given by the
collection of trivial vector bundles  spanned by the vector fields
$x_1\frac\p{\p x_1},\ldots,x_k\frac\p{\p x_k},\frac\p{\p x_{k+1}},\ldots,\frac\p{\p x_n}$
over each chart with the natural transition operators.

For a vector bundle $E$, denote by $\Oc(E)$ the sheaf of its holomorphic sections.
\begin{prop}\label{ehlers}{\em \cite{E}}
There is an exact sequence of coherent sheaves
$$0\to\Oc(V_X(D))\to\Oc(\T X)\to\bigoplus_{i=1}^l\Oc(L_X(D_i))\otimes_{\Oc_X}\Oc_{D_i}\to0.$$

In particular, the tangent bundle $\T X$ has the same Chern classes as the direct sum of the bundle $V_X(D)$
with $L_X(D_1)$,\ldots, $L_X(D_l)$.
\end{prop}

Proposition \ref{ehlers} gives the answer for the Chern classes of $X$, when the Chern classes of $V_X(D)$ are
known. In particular, this is the case when $X$ is a smooth toric
variety, and $D=X\setminus(\C^*)^n$ is the divisor at infinity. In
this case, the vector bundle $V_X(D)$ is trivial, and the Chern
classes of $\T X$ can be found explicitly. This was done by Ehlers
\cite{E}. A more general class of examples is given by regular compactifications
of reductive groups (see the next paragraph for the definition) and, more generally,
of arbitrary spherical homogeneous spaces (see Section \ref{s.sphere}). In
this case, the vector
bundle $V_X(D)$ is no longer trivial but still has a nice
description, which is due to Brion \cite{Brion4}. I recall his result in Subsection \ref{s.chern}
and use it to prove Theorem \ref{compl}.

\paragraph{Regular compactifications.} In this paragraph, I will
define the notion of regular compactifications of reductive groups
following \cite{Brion2}. Let $X$ be a smooth $\DG$--equivariant
compactification of a connected reductive group $G$ of dimension
$n$. Denote by $\Oc_1,\ldots,\Oc_l$ the orbits of codimension one
in $X$. Then the complement $X\setminus G$ to the open orbit is
the union of the closures $\Ocl_1,\ldots,\Ocl_l$ of codimension
one orbits.
\begin{defin}\label{def}
A smooth $\DG$--equivariant compactification $X$ is called regular
if the following three conditions are satisfied.

(1)The hypersurfaces $\Ocl_1,\ldots,\Ocl_l$ are smooth and
intersect each other transversally.

(2)The closure of any $\DG$--orbit in $X\setminus G$ coincides with
the intersection of those hypersurfaces $\Ocl_1,\ldots,\Ocl_l$
that contain it.

(3) For any point $x\in X$ and its $G\times G$--orbit $\Oc_x\subset X$,
the stabilizer $(G\times G)_x\subset \DG$  acts
with a dense orbit on the normal space $T_x X/T_x \Oc_x$ to the orbit.
\end{defin}
This definition was introduced by E.Bifet, De Concini and Procesi in a more general setting
(\cite{BCP}, see also Section \ref{s.sphere}).

If $G$ is a complex torus, then the regularity of $X$ is just
equivalent to the smoothness. However, for other reductive groups,
there exist compactifications that are smooth but not regular. In
particular, it follows from Proposition \ref{Brion} below that
the compactification $X_\pi$ associated with a representation
$\pi:G\to GL(V)$ (see Section \ref{s.cp}) is regular if and only it is smooth and none of the vertices
of the weight polytope of $\pi$ lies on the walls of the Weyl
chambers.

Regular compactifications of reductive groups generalize smooth toric varieties and
retain many nice properties of the latter.
E.g. any regular compactification $X$ can be
covered by affine charts $X_\a\simeq \C^n$ in such a way that only $k$ hypersurfaces
$\Ocl_{i_1}$,\ldots, $\Ocl_{i_k}$ intersect $X_\a$, and
intersections $\Ocl_{i_1}\cap X_\a,\ldots,\Ocl_{i_k}\cap X_\a$ are $k$
coordinate hyperplanes in $X_\a$ \cite{CP2,Brion2}. Here $k$ denotes the rank of $G$. In particular, all
$\DG$-orbits in $X$ have codimension at most $k$, and
all closed orbits have codimension $k$.

If $G$ is of adjoint type, then it has the wonderful
compactification $X_{can}$, which is regular. This example is
crucial for the study of the other regular compactifications.

For arbitrary reductive group $G$, denote by $X_{can}$ the wonderful
compactification of the adjoint group of $G$. There is the
following criterion of regularity.

\begin{prop}{\em \cite{Brion2}}\label{Brion}
Let $X$ be a smooth $\DG$--equivariant compactification of $G$.
Then the condition that $X$ is regular is equivalent to the
existence of a $\DG$--equivariant map from $X$ to $X_{can}$.
\end{prop}
E.g. if $G$ is a complex torus, then the latter condition is
always satisfied because $X_{can}$ is a point in this case.

Thus the set of regular compactifications of $G$ consists of all
smooth $\DG$--equivariant compactifications lying over $X_{can}$. In
particular, for reductive groups of adjoint type the wonderful
compactification is the minimal regular compactification.

\subsection{Demazure bundle and the Chern classes of regular compactifications}\label{s.chern}
In this subsection, I state a formula for the Chern classes of regular compactifications of
reductive groups. It follows from a more general result proved for arbitrary toroidal spherical varieties by
Brion \cite{Brion4}. This formula gives a description of the Chern classes
in terms of two different collections of subvarieties.
The first collection is given by the Chern classes of $G$, which are
independent of a compactification, and the second is given by
the closures of codimension one orbits, which are easy to deal
with (in particular, all their intersection indices with other
divisors  can be computed via the Brion--Kazarnovskii theorem).

Let $X$ be a regular compactification of $G$, and let
$\Ocl_1,\ldots,\Ocl_l$ be the closures of
the $\DG$--orbits of codimension one in $X$. Then the tangent bundle $\T X$ of $X$
can be described using the Demazure vector bundle $V_{can}$ over the wonderful
compactification $X_{can}$ (see Example 1 from Subsection
\ref{s1.chern}) and the line
bundles corresponding to the hypersurfaces $\Ocl_i$.

Let $L(\Ocl_1),\ldots,L(\Ocl_l)$ be the line bundles over $X$
associated with the hypersurfaces $\Ocl_1,\ldots,\Ocl_l$, respectively.
Let $p:X\to X_{can}$ be the canonical map from Proposition \ref{Brion}, and let
$p^*(V_{can})$ be the pull-back of the Demazure vector bundle to $X$. It turns out
that $p^*(V_{can})$ coincides up to a trivial summand
with the logarithmic tangent bundle corresponding to the boundary divisor $X\setminus G$.

\begin{thm}{\em \cite{Brion4}} \label{cor1}
The tangent bundle $\T X$  has the same Chern classes as the direct sum of the
pull-back $p^*(V_{can})$ with the line bundles $L(\Ocl_1),\ldots,L(\Ocl_l)$.
\end{thm}
In the case when $G$ is a complex torus, Theorem \ref{cor1} was proved by Ehlers \cite{E}.
For arbitrary reductive groups,
Theorem \ref{cor1} follows from a more general result by Brion (\cite{Brion4}, 1.6 Corollary~1).

This theorem implies the following formula for the Chern classes $c_1(X),\ldots,c_n(X)$ of
the tangent bundle of $X$. Let $S_i=S_i(TG)\subset G$ for $i=1,\ldots,n-k$ be the Chern classes of the tangent bundle of $G$
defined in the previous section (see Definition \ref{chern}).
Denote by $\o S_i$ the closure of $S_i$ in $X$. Note that $\o S_i$ has proper intersections
with all $\DG$-orbits in $X$ (since this is already true for the wonderful compactification $X_{can}$, and
$X$ lies over $X_{can}$).
\begin{cor} \label{chernclass}
The total Chern class $c(X)=1+c_1(X)+\ldots+c_n(X)$
coincides with the following product:
$$c(X)=(1+\o S_1+\ldots+\o S_{n-k})\cdot\prod_{i=1}^l(1+\Ocl_i).$$
The product in this formula is the intersection product in the (co)homology ring of $X$.
\end{cor}

Below I sketch the proof of Theorem \ref{cor1} following mostly the proofs by Ehlers and Brion.
The goal is to explain the main idea of their proofs, which is very transparent, and motivate
the definition of the Chern classes $S_i$. In the torus case,
this idea can be extended to a complete elementary proof.
For more details see \cite{E} and \cite{Brion4}.

Take $n$ generic vector fields $v_1,\ldots,v_n$  coming from the action of $\DG$.
It is not hard to show that
$v_1,\ldots,v_n$ are generic in the space of all $C^\infty$--smooth
vector fields on $X$ (it is enough to prove it for each affine chart on $X$).
Hence, their degeneracy loci give Chern classes of $X$.
Note that these fields are not only $C^\infty$--smooth but also algebraic
so their degeneracy loci are algebraic subvarieties in $X$.

The picture is especially simple in the torus case,
because in this case $v_1, \ldots, v_{n-i+1}$ are linearly dependent precisely on all orbits
of codimension greater than or equal to $i$ (since they all belong to the tangent bundle of the orbit)
and independent on the other orbits. Hence, the $i$-th Chern class of $X$
consists of all orbits of codimension at least $i$.

In the reductive case, the situation is more complicated because the degeneracy loci of $v_1,\ldots,v_n$ have nontrivial
intersections with the open orbit $G\subset X$. These intersections are exactly the Chern classes
$S_1$,..., $S_{n-k}$ of $G$. So it seems more convenient to use the method described in Subsection \ref{s.prelim}
(see Proposition \ref{ehlers}).
Namely, consider the logarithmic tangent bundle $V_X=V_X(X\setminus G)$  corresponding to the boundary divisor
$X\setminus G=\Ocl_1\cup\ldots\cup\Ocl_l$. Recall that $c$ denotes the dimension of the center of $G$.

\begin{prop}\label{tangent}
The vector bundle $V_X$ is isomorphic to the direct sum of the
pull-back $p^*(V_{can})$ with the trivial vector bundle $E^c$ of
rank $c$.
\end{prop}
\begin{proof}
The vector fields coming from the action of $\DG$ on $X$ are global sections of the bundle $V_X$, since
they are tangent to all $\DG$--orbits in $X$.
It follows easily from condition (3) in the definition of regular compactifications  that these global
sections  span the fiber of $V_X$ at any point of $X$.
Hence, the map $\phi_E:G\to G(n-c,(\gs)/\c)$ considered in Example \ref{Dem2}
extends to a map $p:X\to G(n-c,(\gs)/\c)$.
The rest follows from Example \ref{Dem2}.
\end{proof}

\begin{remark} \em
There is also another construction of the map $p:X\to
X_{can}$ by Brion (see \cite{Brion}).
\end{remark}

\subsection{Applications} \label{proofs}
In this subsection, I prove Theorem
\ref{compl} using the formula for the Chern classes of regular compactifications (Corollary
\ref{chernclass}). Then I apply it to compute the Euler characteristic and the genus of a curve in $G$.

\paragraph{Proof of Theorem \ref{compl}.}
First, define the notion of {\em generic} collection of hyperplane sections
used in the formulation of Theorem \ref{compl}.
 A collection of $m$ hyperplane sections
$H_1$, \ldots, $H_m$ corresponding to representations $\pi_1$,\ldots, $\pi_m$,
respectively,
is called {\em generic}, if
there exists a regular compactification
$X$ of $G$ such that the closure $\o H_i$ of any
hyperplane section $H_i$ is smooth, and all possible intersections
of $\o H_1,\ldots,\o H_m$ with the closures of $\DG$--orbits in $X$ are
transverse.
E.g. one can take the compactification $X_\pi$ corresponding
to the tensor product $\pi$ of the representations $\pi_0,\pi_1,\ldots,\pi_m$, where
$\pi_0$ is any irreducible representation with a strictly dominant highest weight.
Then it is not hard to show that the set of all generic collections
(with respect to the compactification $X_\pi$) is an open dense subset
in the space of all collections.

So the closure $Y=\o C$ of
$C=H_1\cap\ldots\cap H_m$ in $X$ is the transverse intersection of
smooth hypersurfaces. In particular, $Y$ is smooth, and its normal
bundle $N_Y$ in $X$ is the direct sum of $m$ line  bundles
corresponding to the hypersurfaces $\o H_i$. The analogous statement
is true for any subvariety of the form $Y\cap\Ocl_I$, where $I=\{i_1,\ldots,i_p\}$
is a subset of $\{1,\ldots,l\}$ and
$\Ocl_I=\Ocl_{i_1}\cap\ldots\cap\Ocl_{i_p}$.
Let us find the Euler characteristic of
$Y\cap\Ocl_I$ using the classical adjunction formula. Denote by $J=\{1,\ldots,l\}\setminus I$
the complement to the subset $I$. We get that
$\chi(Y\cap\Ocl_I)$ is the term of degree $n$ in the decomposition
of the following intersection product in $X$:

$$(1+\o S_1+\ldots+\o S_{n-k})\cdot\prod_{s=1}^mH_s(1+H_s)^{-1}\cdot\prod
\limits_{i\in I}\Ocl_i\cdot\prod
\limits_{j\in J}(1+\Ocl_j).\eqno(*)$$

On the other hand, since the Euler characteristic is additive, and $C=Y\setminus(\Ocl_1\cup\ldots\cup\Ocl_l)$,
one can express the Euler characteristic
$\chi(C)$ in terms of the Euler characteristics $\chi(Y\cap\Ocl_I)$ over all subsets
$I\subset\{1,\ldots,l\}$:

$$\chi(C)=
\sum_{I\subset\{1,\ldots,l\}}(-1)^{|I|}\chi(Y\cap\Ocl_I).\eqno(**)$$

Combining  formulas $(*)$ and $(**)$, we get the formula of Theorem
\ref{compl}. Indeed, we have that  $\chi(C)$ is the term of degree $n$ in the decomposition
of the following intersection product in $X$:

$$(1+\o S_1+\ldots+\o S_{n-k})\cdot\prod_{s=1}^mH_s(1+H_s)^{-1}\cdot\left(
\sum_{I\sqcup J=\{1,\ldots,l\}}(-1)^{|I|}\prod_{i\in I}\Ocl_i\prod_{j\in J}(1+\Ocl_j)\right).$$

The sum in the parentheses
is equal to $1$, since for any commuting variables $x_1$, $x_2$, \ldots, $x_l$ we have the identity:
$$1=\prod_{i=1}^{l}((1+x_i)-x_i)=\sum_{I\sqcup J=\{1,\ldots,l\}}(-1)^{|I|}\prod_{i\in I}x_i\prod_{j\in J}
(1+x_j).$$

\paragraph{Computation for a curve.}
Apply Theorem \ref{compl} and the formula for the first Chern class $S_1$
to a curve in $G$. We get that if $C=H_1\cap\ldots\cap H_{n-1}$ is a complete intersection of $n-1$
generic hyperplane sections, then
$$\chi(C)=(S_1-H_1-\ldots-H_{n-1})\cdot\prod_{i=1}^{n-1}H_i.$$
Since $S_1$ is also a generic hyperplane section, the computation of $\chi(C)$ reduces to the
computation of the intersection indices of hyperplane sections.

 Recall the Brion-Kazarnovskii
formula for such intersection indices. Denote by $R^+$ the set of all positive roots of $G$. Recall that
$\rho$ denotes the half of the sum of all positive roots of $G$ and $L_T$ denotes
the character lattice of a maximal torus $T\subset G$.
Since $G$ is reductive, we can assume that $\g$ is embedded into $\gl(W)$ so that the trace form ${\rm tr}(A,B)={\rm tr}(AB)$
for $A,B\in\gl(W)$ is nondegenerate on $\g$. Then the inner product $(\cdot ,\cdot)$ on $L_T\otimes\R$ used
in Theorem \ref{degree}
is given by the trace form on $\g$. Choose a Weyl chamber $\D\subset L\otimes\R$.
\begin{thm} {\em \cite{Brion, Kaz}} \label{degree} If $H_\pi$ is a hyperplane section corresponding to a
representation $\pi$ with the
weight polytope $P_\pi\subset L_T\otimes\R$ ,
then the self-intersection index of $H_\pi$ in the ring of conditions is equal to
$$n!\int\limits_{P_\pi\cap\D}\prod_{\a\in R^+}\frac{(x,\a)^2}{(\rho,\a)^2}dx.$$
The  measure $dx$ on $L_T\otimes\R$ is normalized so that the covolume of $L_T$ is $1$.

\end{thm}
This theorem in particular implies that the self-intersection index $H_\pi^n$ depends not on a representation but only
on its weight polytope. Note also that the integrand is invariant under the action of the Weyl group.

Let $H_1, \ldots, H_n$ be $n$ generic hyperplane sections corresponding
to  different representations $\pi_1,\ldots,\pi_n$. To compute their intersection index  one needs
to take the {\em polarization} of $H_\pi^n$. Namely, the formula of Theorem \ref{degree}
gives a polynomial function $D(P)$ of degree $n$ on the space of all virtual polytopes $P\subset L_T\otimes\R$
(the addition in this space is the Minkowski sum). The
{\em polarization} $D_{pol}$ is the unique symmetric
$n$-linear form  on this space such that $D_{pol}(P_\pi,\ldots,P_\pi)=D(P_\pi)$.
Then $D_{pol}(P_{\pi_1},\ldots,P_{\pi_n})$ is the intersection
index $H_1\cdot\ldots\cdot H_n$. For instance, it can be found by applying the differential operator
$\frac1{n!}\frac{\d^n}{\d t_1\ldots\d t_n}$ to the function
$F(t_1,\ldots,t_n)=D(t_1P_{\pi_1}+\ldots+t_nP_{\pi_n})$. E.g. if $P_{\pi_2}=\ldots=P_{\pi_n}$, then
the computation of
$D_{pol}(P_{\pi_1},\ldots,P_{\pi_n})=\left.\frac1{n}\frac{\d}{\d t}\right|_{t=0} D(tP_{\pi_1}+P_{\pi_2})$ reduces to the
integration over the facets of $P_{\pi_2}$.

Thus we get the following answer for $\chi(C)$. For simplicity, the answer is given in the case when
$\pi_1=\ldots=\pi_{n-1}=\pi$. Then its polarization provides the answer in the general case.
Denote by $P_{2\rho}$ the weight polytope of
the irreducible representation of $G$ with the highest weight
$2\rho$.
\begin{cor}\label{curve}
Let $C$ be a curve obtained as the transverse intersection of  a generic collection of $n-1$  hyperplane sections
corresponding to the representation $\pi$.
Then $$\chi(C)=D_{pol}(P_{2\rho},P_\pi,\ldots,P_\pi)-(n-1)D(P_\pi)$$
\end{cor}
A similar answer can be obtained for the genus of $C$ since it is equal to the genus of the compactified curve
$\o C\subset X_{\pi}$. Hence, $g(C)=g(\o C)=1-\chi(\o C)/2$. To compute the Euler characteristic of $\o C$
we need to sum up $\chi(C)$ and the number of points in $\o C\setminus C$. The latter is the intersection index
of $H_\pi^{n-1}$ with the codimension one orbits in $X_\pi$ and can be again computed by the Brion-Kazarnovskii formula.
Choose $l$ facets $F_1,\ldots,
F_l$ of $P_\pi$ so that they parameterize the codimension one orbits in $X_\pi$. This means that each orbit
of the Weyl group acting on the facets of $P_\pi$ contains exactly one $F_i$ (see Theorem \ref{equiv}).
\begin{cor}\label{curve2}
The genus $g(C)$ of $C$ is given by the following formula:
$$g(C)=1-\frac12\left(\chi(C)+(n-1)!\sum_{i=1}^l\int\limits_{F_i\cap\D}\prod_{\a\in R^+}\frac{(x,\a)^2}{(\rho,\a)^2}dx
\right)$$
The measure $dx$ on a facet $F_i$ is normalized as follows. Let $H\subset L\otimes\R$ be the hyperplane containing
$F_i$. Then the covolume of the sublattice $L\cap H$ in $H$ is equal to 1.
\end{cor}

In the above answer, one can rewrite the polarization
$D_{pol}(P_{2\rho},P_\pi,\ldots,P_{\pi})$ in terms of the integrals over
the facets of $P_\pi$.
E.g. in the case when $\pi$ is the irreducible representation with a strictly dominant
highest weight $\lambda$, the answer takes the following form. Let $2\rho=\sum_{i=1}^ka_i\alpha_i$ be the decomposition
of $2\rho$ in the basis of simple roots $\alpha_1,\ldots,\alpha_k$.
$$\chi(C)=n!\left(\frac1{n}\sum_{i=1}^k[a_i\int\limits_{F_i\cap\D}\prod_{\a\in R^+}\frac{(x,\a)^2}{(\rho,\a)^2}dx]-
(n-1)\int\limits_{P_\pi\cap\D}\prod_{\a\in R^+}\frac{(x,\a)^2}{(\rho,\a)^2}dx\right).$$
$$g(C)=1-\frac{n!}2\left(\frac{1}{n}\sum_{i=1}^k[(a_i+1)\int\limits_{F_i\cap\D}\prod_{\a\in R^+}\frac{(x,\a)^2}{(\rho,\a)^2}dx]-
(n-1)\int\limits_{P_\pi\cap\D}\prod_{\a\in R^+}\frac{(x,\a)^2}{(\rho,\a)^2}dx\right)$$

\section{The case of regular spherical varieties} \label{s.sphere}
The results of this paper concerning the Chern
classes of the tangent bundle can be generalized straightforwardly
to the case of arbitrary spherical homogeneous space. In this section, I briefly outline
how this can be done.

Let $G$ be a connected complex reductive group of dimension $r$, and let
$H$ be a closed algebraic subgroup of $G$ . Suppose that the homogeneous space $G/H$ is spherical,
i.e. the action of $G$ on the homogeneous space $G/H$ by left
multiplication is spherical. In the preceding sections, we considered
a particular case of such homogeneous spaces, namely, the space
$(\DG)/G\simeq G$.

The definition of the Chern classes $S_i$ of the tangent bundle $\T(G/H)$
can be repeated
verbatim for $G/H$.
Denote the dimension of $G/H$ by $n$. There is
a space of vector fields on $G/H$ coming from the action of $G$.
Take $n$ arbitrary vector fields $v_1,\ldots,v_n$ of this type.
Define the subvariety $S_i\subset G/H$ as the set of all points
$x\in G/H$ such that the vectors $v_1(x),\ldots,v_{n-i+1}(x)$ are
linearly dependent.

Denote by $\h\subset\g$ the Lie algebra of the subgroup $H$.
Again, there is the Demazure map $p:G/H\to\G(r-n,\g)$, which takes
$g\in G/H$ to the Lie subalgebra $g\h g^{-1}$. Denote by $X_{can}$ the
closure of $p(X)$ in the Grassmannian $\G(r-n,\g)$. This is a
compactification of a spherical homogeneous space $G/N(\h)$, where
$N(\h)\subset G$ is the normalizer of $\h$. Brion conjectured that if $H$ coincides with
$N(H)$, then the compactification $X_{can}$ is smooth, and hence, regular \cite{Brion4}.
F. Knop proved that under the same assumption the normalization of $X_{can}$ is smooth \cite{Knop}.
The conjecture has been proved for semisimple Lie algebras of type A by D.
Luna \cite{Luna}, and  in type D by P. Bravi and G. Pezzini \cite{Bravi}.
In the general case, one can still
define the Demazure bundle over $X_{can}$ as the restriction of
the tautological quotient vector bundle over $\G(r-n,\g)$.

Since we have not used the regularity of $X_{can}$ in the proof of
Lemma \ref{class} the same arguments imply
two facts. First, for a generic choice of vector fields $v_1,\ldots,v_n$,
the resulting subvariety $S_i$ belongs to a fixed class $[S_i]$ in the ring of
conditions. Second, for any compactification $X$ of
$G/H$ lying over $X_{can}$ the closure of a generic $S_i$ in  $X$
intersects properly any orbit of $X$.
Repeating the proof of Lemma \ref{dim} one can also show that $S_i$
is empty unless $i\le n-k$. Here $k$ is the difference between the ranks of $G$ and of $H$.
Therefore, we have $n-k$ well-defined classes $[S_1],\ldots,[S_{n-k}]$ in the ring of
conditions $C^*(G/H)$. Recently, M. Brion and I. Kausz proved that the
{\em $G$--equivariant} Chern classes of the Demazure bundle also vanish for $i>n-k$ \cite{Kausz}.

To extend Theorem \ref{compl} to an arbitrary spherical homogeneous space one can use the same description of the Chern
classes of its regular compactifications. The definition of regular compactifications repeats Definition \ref{def}.
\begin{thm}\label{sphere}
Let $X$ be a regular compactification of $G/H$. Then the total Chern class of $X$ equals to
$$(1+\o S_1+\ldots+\o S_{n-k})\cdot\prod_{i=1}^l(1+\Ocl_i).$$
\end{thm}
This description also follows from Subsection \ref{s.prelim}. The proof uses the methods
mentioned in Subsection \ref{s.chern}. In fact, regular compactifications of spherical homogenous spaces arise naturally
when one try to apply these methods to a wider class of varieties with a group action.
Namely, suppose that a connected complex affine group $G$ of dimension $r$ acts on a compact smooth irreducible complex
variety $X$ with a finite number of orbits.  Then there is a unique open orbit in $X$ isomorphic to $G/H$ for
some subgroup $H\subset G$,
so $X$ can be regarded as a compactification of $G/H$.
Denote by $\Oc_1,\ldots,\Oc_l$
the orbits of codimension one in $X$.
Then one can describe the tangent bundle of $X$ exactly by the methods
mentioned in Subsection \ref{s.chern} if the following conditions
hold. First, the hypersurfaces $\Ocl_1,\ldots,\Ocl_l$ are smooth and
intersect each other transversally (this allows to apply Ehlers' method to the divisor
$X\setminus(G/H)=\Ocl_1\cup\ldots\cup\Ocl_l$). Second, the vector bundle $V_X$ (defined as in Subsection \ref{s.chern})
is generated by its global sections $v_1,\ldots,v_r$, where $v_1,\ldots,v_r$ are infinitesimal
generators of the action of $G$ on $X$ (this allows to give a uniform description of $V_X$ for all
compactifications of $G/H$ satisfying these conditions).
It is not hard to check that these two conditions are equivalent to the definition of regular compactifications.

It turns out that a homogeneous space $G/H$ admits a regular compactification if and only if $G/H$ is spherical
\cite{Bien}.
Regular compactifications of arbitrary spherical homogeneous spaces
are exactly their smooth {\em toroidal} compactifications
\cite{Bien}. A compactification $X$ of the spherical homogeneous space $G/H$ is called {\em toroidal}
if for any codimension one orbit of a Borel subgroup of $G$ acting on $G/H$, its closure in $X$ does not contain
any $G$-orbit in $X$.

The proof of Theorem \ref{compl} goes without any change for
complete intersections in arbitrary spherical homogeneous space
$G/H$. Let $\t H_1,\ldots,\t H_m$ be smooth hypersurfaces in some regular compactification
of $G/H$. Suppose that all possible intersections of $\t H_i$ with the closures of $G$--orbits are transverse.
\begin{thm}\label{compl2}
Let $H_1,\ldots,H_m\subset G/H$ be
the hypersurfaces  $\t H_i\cap(G/H)$, and let $C=H_1\cap\ldots\cap H_m$ be their intersection.
Then the Euler characteristic of $C$ equals to the term of degree $n$ in the decomposition of
$$(1+S_1+\ldots+S_{n-k})\cdot\prod_{i=1}^mH_i(1+H_i)^{-1}.$$
The products are taken in the ring of conditions
$C(G/H)$.
\end{thm}
For instance, if $G/H$ is compact, then the $S_i$ become the usual Chern
classes and the above formula coincides with the classical adjunction formula.
However, if $G/H$ is noncompact then the Chern classes in the
usual sense (as degeneracy loci of {\em generic} vector fields on
$G/H$) do not usually yield the adjunction formula (although they do for $G=(\C^{*})^n$). Indeed, when
the homogeneous space is a noncommutative reductive group, all usual
Chern classes are trivial but as we have seen $\chi(H)\ne(-1)^n
H^n$ even for one smooth hypersurface $H$. Theorem \ref{compl2} shows that
the adjunction formula still holds for noncompact
spherical homogeneous spaces, if one replaces the usual Chern
classes with the refined Chern classes $S_i$ that are defined as
the degeneracy loci of the vector fields coming from the action of
$G$.


\footnotesize
\begin{thebibliography}{99}


\bibitem{Bien}{\sc Fr\'ed\'eric Bien and Michel Brion,} {\em Automorphisms and local rigidity of regular varieties},
Compositio Math. {\bf 104} (1996), no. 1, 1--26
\bibitem{BCP}{\sc E. Bifet, C. De Concini and C. Procesi,} {\em Cohomology of regular embeddings}, Adv. in Math.
{\bf 82} (1990), no. 1, 1--34
\bibitem{Bravi}{\sc P. Bravi and G. Pezzini}, {\em Wonderful varieties of type D}, preprint arXiv.org/math.AG/0410472
\bibitem{Brion}{\sc Michel Brion,} {\em Groupe de Picard et nombres caracteristiques
des varietes spheriques}, Duke Math J. {\bf 58} (1989), no.2,
397--424
\bibitem{Brion4}{\sc Michel Brion,} {\em Vers une generalisation des espaces symetriques}, J. Algebra {\bf 134}
(1990), no. 1, 115--143
\bibitem{Brion2}{\sc Michel Brion,} {\em The behaviour at infinity of the Bruhat decomposition},
Comment. Math. Helv. {\bf 73} (1998), no. 1, 137--174
\bibitem{Kausz}{\sc Michel Brion and Ivan Kausz}, {\em Vanishing of top equivariant Chern classes of regular embeddings}, preprint arxiv.org/math.AG/0503196
\bibitem{C}{\sc C. De Concini,} {Equivariant embeddings of
homogeneous spaces}, Proceedings of the International Congress of
Mathematicians (Berkeley, California, USA, 1986), 369--377
\bibitem{CP2}{\sc C. De Concini and C. Procesi,} {\em Complete symmetric varieties I,}
Lect. Notes in Math. {\bf 996}, Springer, 1983
\bibitem{CP}{\sc  C. De Concini and C. Procesi,} {\em Complete symmetric varieties II Intersection theory,}
Advanced Studies in Pure Mathematics {\bf 6} (1985), Algebraic
groups and related topics, 481--513
\bibitem{E}{\sc F. Ehlers,}
{\it Eine Klasse komplexer Mannigfaltigkeiten und die Aufl\"osung einiger isolierter
Singularit\"aten,} Math. Ann. {\bf 218} (1975), no. 2, 127--157
\bibitem{Fulton2}{\sc W. Fulton,} {\it Intersection
theory,} Springer, 1984
\bibitem{GH}{\sc P. Griffiths, J. Harris,} {\it Principles of algebraic geometry,} Pure and Applied Mathematics,
Wiley-Interscience [John Wiley \& Sons], New York, 1978
\bibitem{GKZ}{\sc I.M. Gelfand, M.M. Kapranov, A.V. Zelevinsky,}
{\em Generalized Euler integrals and A-hypergeometric functions}, Adv. Math. {\bf 84} (1990), no. 2, 255--271
\bibitem{Kapranov2}{\sc M. Kapranov,} {\em Hypergeometric functions on reductive groups,}
Integrable systems and algebraic geometry (Kobe/Kyoto, 1997),
236--281, World Sci. Publishing, River Edge, NJ, 1998
\bibitem{Kiu}{\sc Kiumars Kaveh,} {\em Morse theory and Euler characteristic of sections
of spherical varieties,} Transformation Groups, {\bf 9} (2004), no. 1, 47--63
\bibitem{Kaz}{\sc B.Ya. Kazarnovskii}, {\em Newton polyhedra and the Bezout formula for
matrix-valued functions of finite-dimensional representations,}
{\bf 21} (1987), no. 4, 319--321
\bibitem{K}{\sc Valentina Kiritchenko,} {\em A Gauss-Bonnet theorem, Chern classes and an adjunction formula
for reductive groups}, Thesis, University of Toronto, Toronto, Ontario, 2004
\bibitem{Kleiman}{\sc S.L. Kleiman,} {\em The transversality of a general
translate,} Compositio Mathematica, {\bf 28} (1974), Fasc.3, 287--297
\bibitem{Khov}{\sc A.G. Khovanskii,} {\em Newton polyhedra, and the genus of complete intersections,}
Functional Anal. Appl. {\bf 12} (1978), no. 1, 38--46
\bibitem{Knop2}{\sc  Friedrich Knop,}  {\em The Luna-Vust theory of spherical embeddings,}
Proceedings of the Hyderabad Conference on Algebraic Groups (Hyderabad, 1989), 225--249,
Manoj Prakashan, Madras, 1991
\bibitem{Knop}{\sc Friedrich Knop,} {\em Automorphisms, root systems, and compactifications
of homogeneous varieties,} J. Amer. Math. Soc. {\bf 9} (1996), no. 1, 153--174
\bibitem{Luna} {\sc D. Luna,} {\em Sur les plongements de Demazure,} J. Algebra {\bf 258} (2002),
no. 1, 205--215
\bibitem{Rich} {\sc R.W. Richardson,}  {\it Principal orbit types for algebraic transformation
spaces in characteristic zero,}
Invent. Math. {\bf 16} (1972), 6--14
\bibitem{Rit} {\sc Alvaro Rittatore,} {\it Reductive embeddings are Cohen-Macaulay},
Proc. Amer. Math. Soc. {\bf 131} (2003), no. 3, 675--684
\bibitem{Tim}{\sc D. Timashev,} {\em Equivariant compactifications of reductive groups},
Sb. Math. {\bf 194} (2003), no. 3-4, 589--616
\end{thebibliography}
\end{document}